\documentclass[preprint,12pt]{elsarticle}

\usepackage{amssymb,amsmath,latexsym, amsthm, amscd, mathrsfs, stmaryrd,  color}
\usepackage{epic}
\usepackage{CJK}
\usepackage{setspace}
\usepackage[all]{xy}
\usepackage{yhmath}
\usepackage{soul}
\usepackage{hyperref}
\usepackage{leftidx} 
\usepackage{mathtools} 
\usepackage{rotating} 
\usepackage{stmaryrd} 
\usepackage{youngtab,young}
\usepackage{float}

\usepackage{todonotes}
\usepackage{arydshln}
\makeatletter
\renewcommand*\env@matrix[1][*\c@MaxMatrixCols c]{%
  \hskip -\arraycolsep
  \let\@ifnextchar\new@ifnextchar
  \array{#1}}
\makeatother
\usepackage{pgf,tikz}

\usetikzlibrary{arrows, positioning, calc, chains}
\tikzset{
	ch/.style={circle,draw,on chain,inner sep=2pt},
	chj/.style={ch,join},
	every path/.style={shorten >=4pt,shorten <=4pt}
	}

\numberwithin{equation}{subsection}

\newtheorem{thm}{Theorem}[section]
\newtheorem{cor}[thm]{Corollary}
\newtheorem{lem}[thm]{Lemma}
\newtheorem{rem}[thm]{Remark}

\newtheorem{example}[thm]{Example}

\journal{~~~~~~~~~~}

\begin{document}
\begin{CJK*}{GBK}{song}
\begin{spacing}{1.08}
\begin{frontmatter}

\title{On the Motzkin--Straus theorem and clique number for $k$-uniform hypergraphs }

\author[label1]{Chunli Deng}\ead{chunlideng@163.com}
\author[label1]{Shaohuan Hao}\ead{hsh15075391802@163.com}
\author[label2]{Haifeng Li \corref{cor}}\ead{lihaifeng0820@163.com}

\cortext[cor]{Corresponding author: {Haifeng Li}}

\address{
\address[label1]{Department of Mathematics and Physics, Shijiazhuang Tiedao University, Shijiazhuang 050043, PR China}
\address[label2]{School of Mathematics and Statistics, Guangdong University of Foreign Studies, Guangzhou 510006, PR China}

}

\date{\today}

\begin{abstract}

In this paper, we establish a Motzkin--Straus type equality for  $k$-uniform hypergraphs under certain conditions, which extends the classical Motzkin--Straus theorem for graphs. Using this equality, we present several spectral bounds for the clique number of $k$-uniform hypergraphs in terms of the spectral radii of the adjacency, signless Laplacian, and Laplacian tensors. When $k=2$, the bounds coincide with the corresponding known results for graphs.

\end{abstract}


\end{frontmatter}

\section{Introduction}
A  hypergraph $H=(V(H),E(H))$ consists of a vertex set $V(H)$, which is labeled as $[n]=\{1,\dots,n\}$, and an edge set $E(H)=\{e_1,\dots,e_m\}$ with
$e_p\subseteq V(H)$ for $p\in[m]$.
 $H$ is \textit{$k$-uniform} if each   $e \in E(H)$ contains exactly $k$ vertices.
In particular, a $2$-uniform hypergraph is called a graph.
 For a $k$-uniform hypergraph $H$, a subset $S \subseteq V(H)$ is called a \textit{$t$-clique} if every $t$ vertices in $S$ are contained in some edge of
$H$, $t=1,\ldots,k$. The largest cardinality of a $t$-clique of $H$ is called the \textit{$t$-clique number} of $H$, denoted by $\omega_t(H)$, $t=1,\ldots,k$. The  \textit{$k$-clique number} $\omega_k(H)$ is also called the clique number of $H$, and is simply denoted  $\omega$.

Let $\mathfrak{S}=\left\{\mathbf{x}=(x_1,x_2,\dots,x_n)^\top\in \mathbb{R}_+^n\,~\bigg|\,~\sum\limits_{i=1}^n x_i=1\right\}$, where $\mathbb{R}_+^n$ is the set of $n$-dimensional nonnegative vectors.
 The \textit{Lagrangian of a graph} $G=(V(G),E(G))$ over the $\mathfrak{S}$ denoted by
\[
L(G,\mathbf{x})=\sum_{\{i,j\}\in E(G)}x_ix_j.
\]
Let $L(G)$ be
 the largest Lagrangian of $G$, i.e.,
\[
L(G) = \max_{\mathbf{x}\in \mathfrak{S}}L(G,\mathbf{x}).
\]

In 1965, Motzkin and Straus \cite{1} established a remarkable connection between the Lagrangian  of a graph and its clique number as follows.
\begin{thm} \cite{1}\label{th1.1}
Let $G$ be a graph with clique number $\omega$. Then
\[
L(G) = \frac12\left(1-\frac{1}{\omega}\right).
\]
\end{thm}

This result also provided a new proof of Tur\'an's theorem~\cite{1} and has had an important  impact on the study of clique-related problems.
In spectral graph theory, the Motzkin--Straus theorem and its extensions have been used to derive spectral bounds for the clique number [2-5].
It has also been applied to the maximum clique problem in combinatorial optimization~\cite{Pelillo}. Related results on Tur\'an-type problems in spectral and extremal graph and hypergraph   theory have also been obtained~\cite{Keevash,liu-ning}.

In \cite{10}, Wilf presented a lower bound for the  clique
number $\omega$ in terms of the spectral radius $\rho(A)$ of the adjacency matrix $A$ of the graph $G$.
\begin{thm}\cite{10}\label{th1.2}
Let $G$ be a graph with clique number $\omega$. Then
\[
\omega \ge \frac{\|\mathbf{u}\|_1^2}{\|\mathbf{u}\|_1^2 - \rho(A)},
\]
where $\mathbf{u}$ is the nonnegative eigenvector corresponding to $\rho(A)$ with $\|\mathbf{u}\|_2=1$.

\end{thm}


In \cite{3}, Lu et al presented another lower bound for the  clique
number $\omega$ by the spectral radius $\rho(L)$ of the Laplacian matrix $L$ of the graph $G$.
\begin{thm} \cite{3}\label{th1.3}
Let $G$ be a graph with clique number $\omega$, $m$ edges and maximum degree $\Delta$. Then
\[
\omega \ge \frac{2m}{2m-(\rho(L)-\Delta)^2}.
\]
\end{thm}

Motivated by these results, we want to study the link between the Lagrangian of hypergraph  and clique number and give the applications in spectral hypergraph theory. The \textit{Lagrangian of a hypergraph} $H=(V(H),E(H))$ over the $\mathfrak{S}$ denoted by
\[
L(H,\mathbf{x}) = \sum_{\substack{ \{i_1,\ldots,i_k\} \in E(H)}} x_{i_1}\cdots x_{i_k},
\]
and $L(H)$ is the largest Lagrangian of $H$, i.e.,
\[
L(H) = \max_{\mathbf{x}\in \mathfrak{S}}L(H,\mathbf{x}).
\]

However, the Motzkin--Straus theorem does not extend directly to
uniform hypergraphs in general.
Consequently, a lot of effort has been devoted to developing
Motzkin--Straus type results and studying Lagrangians of hypergraphs.
Frankl et al.~\cite{13} formulated a conjecture on the maximum Lagrangian of an $r$-uniform hypergraph with a prescribed number of edges, which was verified in certain cases by Talbot~\cite{5} and further developed by Peng et al.~\cite{9}. Peng et al.~\cite{12} also extended the Motzkin--Straus theorem to a class of non-uniform hypergraphs.

Spectral methods have also been used to study clique numbers of
hypergraphs.
 Xie et al.\cite{26}
gave inequalities relating the Laplacian and signless Laplacian
$H$-eigenvalues to the clique and coclique numbers of uniform
hypergraphs.
Hou et al.~\cite{27} studied homogeneous polynomials and clique-related problems for general hypergraphs.
Duan et al. \cite{25} established bounds on the
spectral radius of general hypergraphs in terms of the clique number.
Further studies on spectral properties of hypergraphs and related spectral methods can be found in [16-21].

Motivated by these studies, we first establish a Motzkin--Straus type equality for  $k$-uniform hypergraphs under certain conditions. Based on this equality, we further obtain  corresponding bounds for the clique number in terms of the spectral radii of the adjacency, signless Laplacian, and Laplacian tensors. When $k=2$, these results coincide with Theorems~\ref{th1.1}, \ref{th1.2}, and~\ref{th1.3} for graphs. Hence, our results extend the corresponding results for graphs to $k$-uniform hypergraphs.

\section{Preliminaries}

In this section, we introduce the basic definitions and preliminary results needed in the sequel.

A $k$-th order $n$-dimensional tensor $\mathcal{T}=(t_{i_1{i_2} \cdots {i_k}}) ~({i_j\in[n],~j\in[k]})$ is a multidimensional array with $n^k$ entries.
If $t_{i_1{i_2} \cdots {i_k}}\geq 0$ for all ${i_j\in[n],~j\in[k]}$, then $\mathcal{T}$ is called \textit{nonnegative}.
If  $t_{i_{1}i_{2}\cdots i_{k}}=t_{\sigma(i_{1})\sigma(i_{2})\cdots \sigma(i_{k})}$, where $\sigma$ is any permutation of the indices $i_1,\ldots,i_k$, then $\mathcal{T}$ is called \textit{symmetric}.
For a $k$-uniform hypergraph $H$ with $n$ vertices,
the \textit{adjacency tensor} \cite{14} $\mathcal{A}=(a_{i_1 i_2  \cdots i_k }) ~(i_1,\ldots , i_k\in[n])$ of $H$ is a  $k$-th order $n$-dimensional nonnegative symmetric tensor, and the entries
\[
a_{i_1 i_2  \cdots i_k } =
          \begin{cases}\label{eq1.2}
            \frac{1}{(k-1)!}, &{\text{if} ~\{i_1 ,i_2 , \ldots ,i_k \}\in E({{H}}),} \\
            0,  &{\text{otherwise}.}
          \end{cases}
\]
 For a vertex $i\in V(H)$, the \textit{degree}  $d_i=|E_i|$, where  $E_i=\{e\in E(H)\mid i\in e\}$.
Let $\mathcal{D}$ be a $k$-th order $n$-dimensional diagonal tensor with its diagonal element $d_{i\cdots i}$ being $d_i$, the degree
of vertex $i$ in $H$, for all $i\in[n]$.
Then $\mathcal{L}=\mathcal{D}-\mathcal{A}$ is the \textit{Laplacian tensor}, and $\mathcal{Q}=\mathcal{D}+\mathcal{A}$ is
the \textit{signless Laplacian tensor} of $H$ \cite{6}.

Given a  $k$-th order $n$-dimensional complex tensor $\mathcal{T}=(t_{i_1\cdots i_k})$
and a complex vector $\mathbf{x}=(x_1,x_2,\dots,x_n)^\top$.
Denote $\mathcal{T}\mathbf{x}^{k-1}$   a vector in $\mathbb{C}^n$,  its $i$th component defined by
\[
\bigl(\mathcal{T}\mathbf{x}^{k-1}\bigr)_i = \sum_{i_2,\dots,i_k=1}^n t_{ii_2\cdots i_k}x_{i_2}\cdots x_{i_k}.
\]
If there exist $\mathbf{x}\in\mathbb{C}^n$ ($\mathbf{x}\neq 0$) and $\lambda\in\mathbb{C}$ satisfy
\[
\mathcal{T}\mathbf{x}^{k-1}=\lambda \mathbf{x}^{[k-1]},
\]
then  $\lambda$ is called an \textit{eigenvalue} of $\mathcal{T}$, and $\mathbf{x}$ is called its corresponding \textit{eigenvector} of $\mathcal{T}$ \cite{Qi05}, where $\mathbf{x}^{[k-1]}=(x_1^{k-1},\ldots,x_n^{k-1})^\top$. We use $\rho(\mathcal{T})$ to represent the spectral radius of $\mathcal{T}$, which is given by the largest modulus among the eigenvalues of $\mathcal{T}$.

Perron--Frobenius theory plays a fundamental role in the study of nonnegative matrices and tensors.
\begin{lem}\cite{yang-perron}\label{perron-Frobenius}
If $\mathcal{T}$ is a  $k$-th order $n$-dimensional nonnegative tensor, then $\rho(\mathcal{T})$ is a nonnegative eigenvalue with a nonnegative eigenvector
 corresponding to it.
\end{lem}

In this paper, we denote the spectral radii of the
adjacency, Laplacian, and signless Laplacian tensors of the $k$-uniform hypergraph $H$ by $\rho(\mathcal{A})$, $\rho(\mathcal{L})$, $\rho(\mathcal{Q})$, respectively.

\begin{lem}\cite{6}\label{Qxiajie}
 Let $H$ be a $k$-uniform hypergraph with maximum degree  $\Delta>0$ and average degree $\bar{d}$. Then
\[
\max\left\{\Delta,2\bar{d}\right\} \le \rho(\mathcal{Q}) .
\]
\end{lem}

\begin{lem}\cite{8}\label{QandL}
Let $H$ be a $k$-uniform hypergraph.  Then
\[
\rho(\mathcal{L}) \le \rho(\mathcal{Q}).
\]
\end{lem}

\begin{lem}\cite{6} \label{Qi}
Let $H$ be a $k$-uniform hypergraph with  maximum degree $\Delta>0$. Then
\[
\rho(\mathcal{L})\ge \Delta.
\]
\end{lem}

The complete $k$-uniform hypergraph on $n$ vertices, denoted by
$K_n^{k}$, is a $k$-uniform hypergraph whose edge set  consists of all $k$-element subsets of the vertex set.
Obviously, for  complete $k$-uniform hypergraph $K_n^k$, its clique number equals the vertex number $n$. By Maclaurin's inequality,
\begin{align}\label{complete}
L(K_n^k,\mathbf{x})=\sum\limits_{{ \{i_1,\ldots,i_k\} \in E(K_n^k)}} x_{i_1}\cdots x_{i_k}
\leq \left(\frac{\sum\limits_{i=1}^n{x_i}}{n}\right)^k \binom{n}{k}
=\left(\frac{1}{n}\right)^k \binom{n}{k}, \tag{2.1}
\end{align}
and equality holds if and only if $x_1=\cdots =x_n=\frac{1}{n}$.
Thus $L(K_n^k,\mathbf{x})$ attains its maximum when $x_1=\cdots = x_n=\frac{1}{n}$. Consequently, we have
\[
L(K_n^k)=\left(\frac{1}{n}\right)^k\binom{n}{k}.
\]

In \cite{9}, Peng et al. gave an extension of the Motzkin--Straus theorem  to hypergraphs under certain conditions.
\begin{lem}\cite{9}\label{peng}
Let $H$ be a $k$-uniform hypergraph with $m$ edges and clique number $\omega$.
If $m \le \binom{\omega}{k} + 2(\omega+1 - k)$, then
\[
L(H) = L(K_\omega^k)=\left(\frac{1}{\omega}\right)^k \binom{\omega}{k}.
\]
\end{lem}

For any vector $\mathbf{x}=(x_1,\dots,x_n)^\top\in\mathbb{R}^n$, the \textit{support} of $\mathbf{x}$, denoted by
$\operatorname{supp}(\mathbf{x})$, is defined as
$\operatorname{supp}(\mathbf{x}) = \{i\in [n] \mid x_i \neq 0\}$.

Throughout this paper, all $k$-uniform hypergraphs are assumed to be nonempty and $k\geq 2$. Hence, their $t$-clique numbers satisfy $\omega_t(H)\geq k$, $t=2,\ldots,k$.

\section{A Motzkin--Straus type result for uniform hypergraphs}

The classical Motzkin--Straus theorem establishes a fundamental connection
between the clique number and the Lagrangian of a graph. In this section, we extend this
connection to  uniform hypergraphs under certain conditions.  We then present two examples
to illustrate the role of the condition in the theorem and to compare it
with an existing related condition.

\begin{thm}\label{th3.1}
Let $H$ be a $k$-uniform hypergraph with clique number $\omega$. If $\omega_2(H)=\omega$,
then
\begin{align*}
    L(H)=\left(\frac{1}{\omega}\right)^k \binom{\omega}{k}.
\end{align*}
\end{thm}
\begin{proof}
Assume $H_0$ is the maximum complete $k$-uniform subhypergraph of $H$ on $\omega$ vertices.
Let $\mathbf{y} = (y_1, \dots, y_n)^\top \in \mathbb{R}_+^n$  with entries
\[
y_i =
\begin{cases}
\dfrac{1}{\omega} ,& \text{if } i \in V(H_0), \\[6pt]
0 ,& \text{otherwise}.
\end{cases}
\]
It is easy to verify that $\sum\limits_{i=1}^{n} y_i = 1$.
Then we have
\begin{align*}
L(H) &\ge L(H,\mathbf{y}) \\
&= \sum_{e\in E(H_0)}y_{i_1}y_{i_2} \cdots y_{i_k} \\
&= \left(\frac{1}{\omega}\right)^k \binom{\omega}{k}.
\end{align*}

Next, we prove $L(H) \le \left(\frac{1}{\omega}\right)^k \binom{\omega}{k}$.

Assume that $\mathbf{x}^*= (x_1, \dots, x_n)^\top \in \mathbb{R}_+^n$ is a vector that maximizes $L(H)$ with the smallest support set, and let $S$ be the support set
corresponding to $\mathbf{x}^*$.

First, we prove that $S$ is a $2$-clique by contradiction. Assume that $u, v \in S$, and $\{u, v\} \not\subseteq e$ for any $e \in
E(H)$. Then
\[
L(H) = L(H, \mathbf{x}^*) = A x_u + B x_v + C,
\]
where $A=\sum\limits_{\{u,i_2,\dots,i_k\} \in E_u} x_{i_2} \cdots x_{i_k}, B=\sum\limits_{\{v,i_2,\dots,i_k\} \in E_v} x_{i_2} \cdots
x_{i_k}, C = \sum\limits_{\substack{\{i_1,\dots,i_k\} \in E(H), \\ u, v \notin \{i_1,\dots,i_k\}}} x_{i_1} \cdots x_{i_k}.$

If $A \ge B$, let $\mathbf{z} = (z_1, \dots, z_n)^\top \in \mathbb{R}_+^n$, where
\[
z_i =
\begin{cases}
x_u + x_v, & i = u, \\[4pt]
0, & i = v, \\[4pt]
x_i, & i \neq u, v.
\end{cases}
\]
Then $L(H, \mathbf{z}) = A (x_u +  x_v) + C$.
Thus $L(H) =L(H, \mathbf{x}^*) \le L(H, \mathbf{z}) \le L(H)$. It yields that $\mathbf{z}$ is also a vector that maximizes $L(H)$, and $|\operatorname{supp}(\mathbf{z})| < |\operatorname{supp}(\mathbf{x}^*)|$,  which contradicts the definition of $\mathbf{x}^*$. Similarly, a contradiction
is also obtained when $A < B$.
Therefore, $S$ is a $2$-clique. Clearly, $|S| \geq k$.

Let the $|S|$-dimensional vector $\boldsymbol{\mathbf{x}}'$ be obtained from the $n$-dimensional vector $\boldsymbol{\mathbf{x}^*}$ by
removing its zero components. Then  we obtain
\begin{align}\label{eq3.1}
L(H) = L(H, \mathbf{x}^*) \leq L(K_{|S|}^{k},\mathbf{x}^\prime) \le L\left(K_{|S|}^{k}\right).\tag{3.1}
\end{align}
It follows from (\ref{complete}) that
\begin{align}\label{eq3.2}
L\left(K_{|S|}^{k}\right) =\left(\frac{1}{|S|}\right)^k\binom{|S|}{k}.\tag{3.2}
\end{align}

Since $S$ is a $2$-clique, we get $|S|\leq\omega_2(H) =\omega$. Thus
\begin{align}\label{eq3.3}
\left(\frac{1}{|S|}\right)^k\binom{|S|}{k}&=\frac{1}{|S|^k}\left(\frac{|S|(|S|-1)\cdots(|S|-k+1)}{k!} \right) \notag\\
&=\frac{1}{k!}(1-\frac{1}{|S|})(1-\frac{2}{|S|})\cdots (1-\frac{k-1}{|S|}) \notag\\
&\leq\frac{1}{k!}(1-\frac{1}{\omega})(1-\frac{2}{\omega})\cdots (1-\frac{k-1}{\omega}) \notag\\
&=\frac{1}{\omega^k}\left(\frac{\omega(\omega-1)\cdots(\omega-k+1)}{k!} \right)\notag\\
&=\left(\frac{1}{\omega}\right)^k \binom{\omega}{k}.\tag{3.3}
\end{align}
Then equations (\ref{eq3.1})-(\ref{eq3.3}) yield that
\[
L(H) \le \left(\frac{1}{\omega}\right)^k \binom{\omega}{k}.
\]

In conclusion,
\[
L(H) = \left(\frac{1}{\omega}\right)^k \binom{\omega}{k}.
\]

\end{proof}

\begin{rem}
In Theorem \ref{th3.1}, when $k=2$, a $2$-uniform hypergraph is   a graph,
and a $2$-clique is exactly a clique. Hence
\[
\omega_2(G)=\omega(G)
\]
holds automatically for every graph $G$. Therefore, Theorem
\ref{th3.1} yields
\[
L(G)
=
\left(\frac{1}{\omega}\right)^2 \binom{\omega}{2}
=
\frac{1}{2}\left(1-\frac{1}{\omega}\right),
\]
which is precisely the classical Motzkin--Straus theorem, i.e. Theorem \ref{th1.1}. Thus, Theorem
\ref{th3.1} can be regarded as a Motzkin--Straus type extension from
graphs to $k$-uniform hypergraphs.
\end{rem}

Next, we give an example to show that the equality in Theorem~\ref{th3.1} need not hold when the condition $\omega_2(H)=\omega$ is not satisfied.
\begin{example}
Let $H$ be the $3$-uniform hypergraph with vertex set
$
V(H)=\{1,2,3,4,5\}
$
and edge set
\[
E(H)=
\big\{
\{1,3,4\},\{1,3,5\},\{1,4,5\},
\{2,3,4\},\{2,3,5\},\{2,4,5\}
\big\}.
\]
It is easy to verify that the clique number of $H$ is
$
\omega(H)=3
$
and
$
\omega_2(H)=4.
$
So $H$ does not satisfy the condition of Theorem \ref{th3.1}.

For any  vector
$
\mathbf{x}=(x_1,x_2,x_3,x_4,x_5)^{\top}$ with $
x_i\geq 0$ and $\sum\limits_{i=1}^5x_i=1,
$
we have
\begin{align*}
L(H,\mathbf{x})
&=
x_1x_3x_4+x_1x_3x_5+x_1x_4x_5
+x_2x_3x_4+x_2x_3x_5+x_2x_4x_5 \\
&=
(x_1+x_2)(x_3x_4+x_3x_5+x_4x_5).
\end{align*}
Set
$
a=x_1+x_2.
$
Then
$
x_3+x_4+x_5=1-a.
$
By
$
x_3x_4+x_3x_5+x_4x_5
\leq
\frac{(x_3+x_4+x_5)^2}{3},
$
we obtain
\[
L(H,\mathbf{x})
\leq
\frac{a(1-a)^2}{3}.
\]
Since
$
\frac{a(1-a)^2}{3}
~(0\leq a\leq 1)
$
attains its maximum at $a=\frac13$, it follows that
\[
L(H)\leq \frac{4}{81}.
\]
Equality holds when
$
x_1+x_2=\frac13
$
and
$
x_3=x_4=x_5=\frac29.
$
For instance, taking
$
\mathbf{x}
=
\left(
\frac16,\frac16,\frac29,\frac29,\frac29
\right)^{\top}
$
gives
\[
L(H)=\frac{4}{81}.
\]
Therefore,
\[
L(H)
=
\frac{4}{81}
>
\frac{1}{27}
=
\left(\frac{1}{\omega}\right)^3
\binom{\omega}{3}.
\]
Thus, the conclusion of Theorem \ref{th3.1} does not hold for this
hypergraph.

\vspace{3mm}
Now add the edge $\{3,4,5\}$ to $H$, and denote the resulting
$3$-uniform hypergraph by $H_0$. Thus,
$
V(H_0)=\{1,2,3,4,5\}
$
and
\[
E(H_0)=
\big\{
\{1,3,4\},\{1,3,5\},\{1,4,5\},
\{2,3,4\},\{2,3,5\},\{2,4,5\},
\{3,4,5\}
\big\}.
\]
It is easy to check that
$
\omega(H_0)=4,
$
and that $H_0$ satisfies the condition of Theorem \ref{th3.1}.
Then we get
\[
\left(\frac{1}{\omega}\right)^3
\binom{\omega}{3}
=
\left(\frac14\right)^3
\binom43
=
\frac1{16}.
\]

For any  vector
$
\mathbf{x}=(x_1,x_2,x_3,x_4,x_5)^{\top}$ with $
x_i\geq 0$ and $\sum\limits_{i=1}^5x_i=1,
$
\begin{align*}
L(H_0,\mathbf{x})
&=
(x_1+x_2)(x_3x_4+x_3x_5+x_4x_5)
+x_3x_4x_5.
\end{align*}
Again, let
$
a=x_1+x_2,
$
so that
$
x_3+x_4+x_5=1-a.
$
Using
$
x_3x_4+x_3x_5+x_4x_5
\leq
\frac{(1-a)^2}{3}
$
and
$
x_3x_4x_5
\leq
\left(\frac{1-a}{3}\right)^3
=
\frac{(1-a)^3}{27},
$
we obtain
\begin{align*}
L(H_0,\mathbf{x})
&\leq
\frac{a(1-a)^2}{3}
+
\frac{(1-a)^3}{27} =
\frac{(1-a)^2(8a+1)}{27}.
\end{align*}
By calculation,  $\frac{(1-a)^2(8a+1)}{27}$ attains its maximum at
$
a=\frac{1}{4},
$
and consequently
\[
L(H_0)\leq \frac1{16}.
\]
Equality holds when
$
x_1+x_2=\frac{1}{4}
$
and
$
x_3=x_4=x_5=\frac{1}{4}.
$
For instance, taking
$
\mathbf{x}
=
\left(
\frac18,\frac18,\frac14,\frac14,\frac14
\right)^{\top}
$
gives
\[
L(H_0)=\frac1{16}.
\]
Therefore,
\[
L(H_0)
=
\frac1{16}
=
\left(\frac{1}{\omega}\right)^3
\binom{\omega}{3}.
\]

This shows that for $k$-uniform hypergraph $H$ satisfying the condition of Theorem \ref{th3.1}, 
the calculation result of $L(H)$ is equal to $\left(\frac{1}{\omega}\right)^k
\binom{\omega}{k}$.

\end{example}

In the following example, we compare the conditions of
Lemma \ref{peng} and Theorem \ref{th3.1} and show that neither
of them implies the other.

\begin{example}
Let $H_1$ be the $3$-uniform hypergraph with $9$ vertices and $4$ edges
shown in Figure \ref{fig:H_example}.

\begin{figure}[H]
\centering
\includegraphics[width=0.6\textwidth]{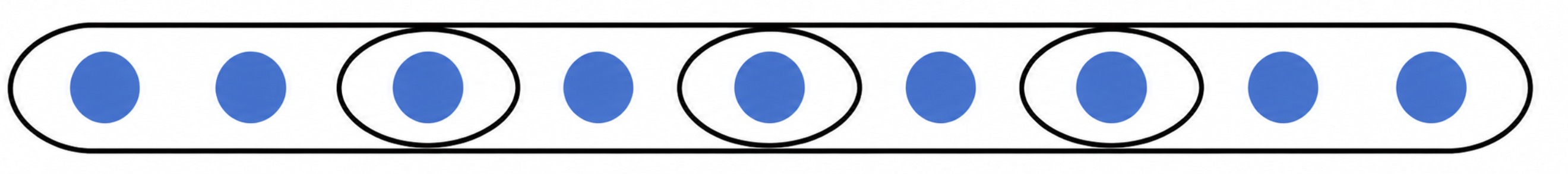}
\caption{A $3$-uniform hypergraph $H_1$.}
\label{fig:H_example}
\end{figure}

It is straightforward to verify that $H_1$ satisfies the condition of
Theorem \ref{th3.1}, but does not satisfy the condition of
Lemma \ref{peng}.

On the other hand, let $H_2$ be the $3$-uniform hypergraph with $4$ vertices and $3$ edges
shown in Figure  \ref{fig:P_example}.

\begin{figure}[H]
\centering
\includegraphics[width=0.25\textwidth]{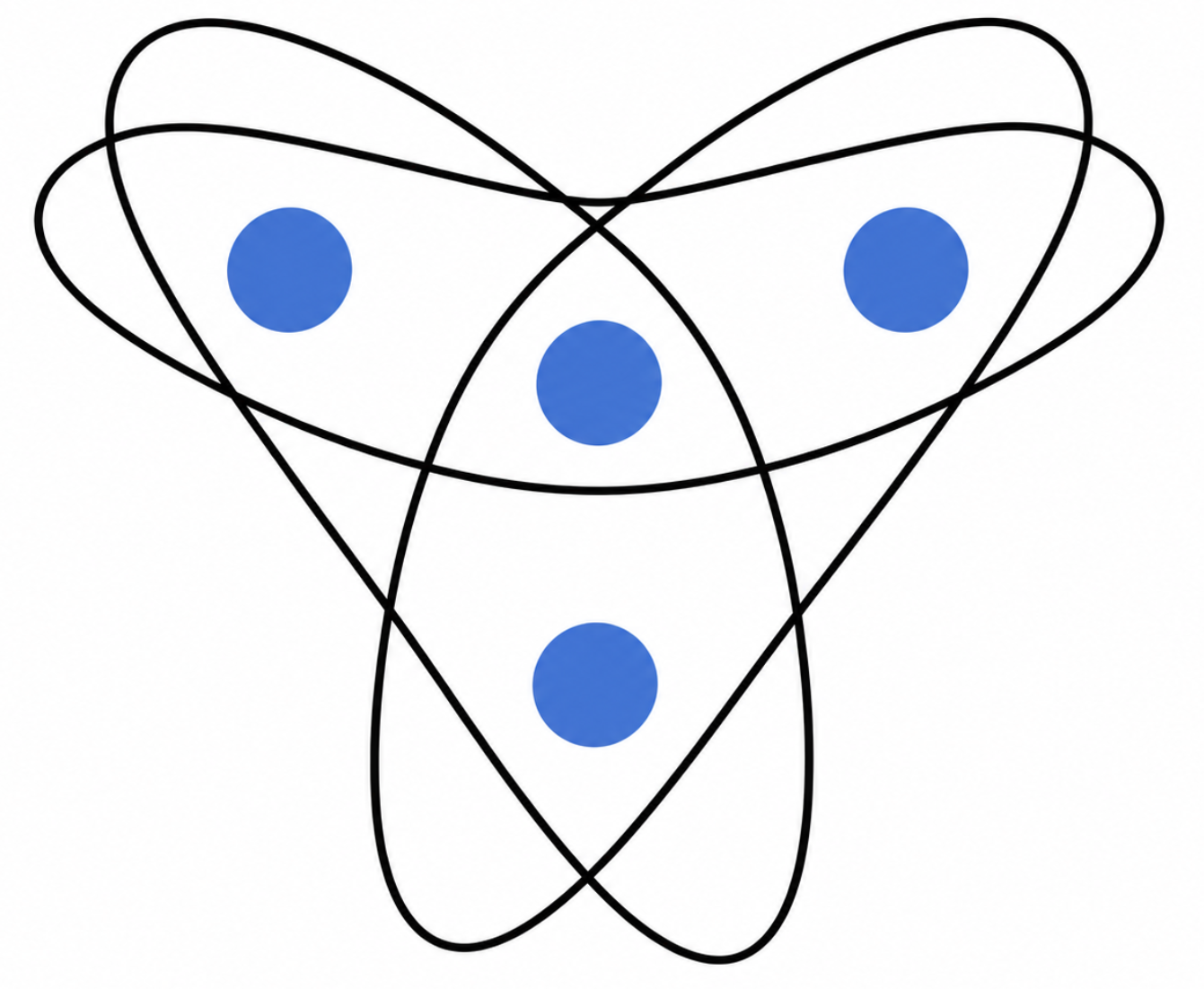}
\caption{A $3$-uniform hypergraph $H_2$.}
\label{fig:P_example}
\end{figure}

It is easy to check that $H_2$ satisfies the condition of
Lemma \ref{peng}, but does not satisfy the condition of
Theorem \ref{th3.1}.

Therefore, the conditions of Lemma \ref{peng} and
Theorem \ref{th3.1} are incomparable; that is, neither condition
implies the other.
\end{example}

\section{Spectral bounds for the clique number of uniform hypergraphs }

In this section, we establish several  lower bounds for the clique number of a $k$-uniform hypergraph in terms of the spectral radii of adjacency, signless Laplacian, and Laplacian tensors.

\subsection{Bounds based on the adjacency tensor}

By Lemma \ref{perron-Frobenius}, the spectral radius
$\rho(\mathcal{A})$ of the adjacency tensor $\mathcal{A}$ is nonnegative and  there is  a nonnegative eigenvector corresponding to it. Let $\mathbf{u}=(u_1,u_2,\dots,u_n)^\top $ be a nonnegative eigenvector corresponding to $\rho(\mathcal{A})$ and
  $\|\mathbf{u}\|_k= (\sum\limits_{i=1}^{n} u_i^k)^{\frac{1}{k}}=1 $.
Then we have the following conclusion, which extends the Theorem \ref{th1.2} from graphs to $k$-uniform hypergraphs.
\begin{thm}\label{th4.1}
Let $H$ be a $k$-uniform hypergraph with clique number $\omega$. If $\omega_2(H)=\omega$,
then
\[
   \omega \ge \frac{1}{1 - \sqrt[k-1]{\dfrac{(k-1)!\rho(\mathcal{A})}{\|\mathbf{u}\|_1^k}}},
\]
where $\mathbf{u}$ is the nonnegative eigenvector corresponding to $\rho(\mathcal{A})$ with $\|\mathbf{u}\|_k=1$.
\end{thm}
\begin{proof}
By Theorem \ref{th3.1}, $\max\limits_{\mathbf{x}\in \mathfrak{S}}\sum\limits_{\substack{ \{i_1,\ldots,i_k\} \in E(H)}} x_{i_1}\cdots x_{i_k}=\left(\frac{1}{\omega}\right)^k \binom{\omega}{k}.$
Let $\mathbf{y}=\frac{\textbf{u}}{\|\mathbf{u\|_1}}$, then  $\mathbf{y} \in \mathfrak{S}$. Thus, we have
\begin{align*}\label{eq4.1}
k\left(\frac{1}{\omega}\right)^k \binom{\omega}{k}&\ge k\sum_{\{i_1,i_2,\dots,i_k\}\in E(H)} y_{i_1}y_{i_2}\cdots y_{i_k}
 \\
 &= \sum_{i_1,i_2,\dots,i_k\in [n]} a_{i_1i_2\dots i_k}y_{i_1}y_{i_2}\cdots y_{i_k}
 \\
&= \sum_{i_1,i_2,\dots,i_k\in [n]} a_{i_1i_2\dots i_k}\frac{u_{i_1}}{\|\mathbf{u}\|_1}\frac{u_{i_2}}{\|\mathbf{u}\|_1}\cdots
\frac{u_{i_k}}{\|\mathbf{u}\|_1} \\
&= \frac{1}{\|\mathbf{u}\|_1^k} \mathbf{u}^\top \mathcal{A}\mathbf{u}^{k-1}\\
&= \frac{1}{\|\mathbf{u}\|_1^k} \mathbf{u}^\top \rho(\mathcal{A})\mathbf{u}^{[k-1]}\\
&= \frac{\rho(\mathcal{A})}{\|\mathbf{u}\|_1^k} \mathbf{u}^\top \mathbf{u}^{[k-1]}\\
&= \frac{\rho(\mathcal{A})}{\|\mathbf{u}\|_1^k}. \tag{4.1}
\end{align*}
And
\begin{align*}\label{eq4.2}
\left( \frac{1}{\omega} \right)^k \binom{\omega}{k}
&= \frac{1}{\omega^k} \cdot \frac{\omega(\omega-1)\cdots(\omega-k+1)}{k!} \\
&= \frac{1}{k!} \cdot \frac{\omega}{\omega} \cdot \frac{\omega-1}{\omega} \cdots \frac{\omega-k+1}{\omega} \\
&\le \frac{1}{k!} \cdot \left( 1 - \frac{1}{\omega} \right)^{k-1}.\tag{4.2}
\end{align*}
Then from inequalities (\ref{eq4.1}) and (\ref{eq4.2}), we get
\begin{align*}
   \omega \ge \frac{1}{1 - \sqrt[k-1]{\dfrac{(k-1)!\rho(\mathcal{A})}{\|\mathbf{u}\|_1^k}}}.
\end{align*}

\end{proof}

\begin{rem}
When $k=2$, Theorem  \ref{th4.1} coincides with  Theorem \ref{th1.2} for graphs. Therefore, Theorem \ref{th4.1} extends the result of graphs  to  $k$-uniform hypergraphs  satisfying $\omega_2(H)=\omega$.
\end{rem}

Next, we give another lower bound for the clique number of a
$k$-uniform hypergraph in terms of the spectral radius of its
adjacency tensor $\rho(\mathcal{A})$.

\begin{thm}\label{th4.3}
Let $H$ be a $k$-uniform  hypergraph with $n$ vertices, $m$ edges and  clique number $\omega$. If $\omega_2(H)=\omega$,
then
\[
\omega \ge \frac{1}{1-\sqrt[k-1]{k!m^{1-k}\left(\frac{\rho(\mathcal{A})}{k}\right)^k}}.
\]
\end{thm}
\begin{proof}
Let $\mathbf{u}=(u_1,u_2,\dots,u_n)^\top $ be a nonnegative eigenvector corresponding to $\rho(\mathcal{A})$ and
  $ \sum\limits_{i=1}^{n} u_i^k=1 $. So
\begin{align*}
\rho(\mathcal{A}) = \mathbf{u}^\top \mathcal{A} \mathbf{u}^{k-1} = k\sum_{\{i_1,i_2,\dots,i_k\}\in E(H)} u_{i_1}u_{i_2}\cdots u_{i_k}.
\end{align*}
By the Power Mean inequality, we have
\begin{align*}
\frac{\sum\limits_{\{i_1,i_2,\dots,i_k\}\in E(H)} u_{i_1}u_{i_2}\cdots u_{i_k}}{m}
\le
\left(
\frac{\sum\limits_{\{i_1,i_2,\dots,i_k\}\in E(H)} u_{i_1}^k u_{i_2}^k \cdots u_{i_k}^k}{m}
\right)^{\frac{1}{k}}.
\end{align*}
Thus,
\[
\rho(\mathcal{A}) \le k m^{1-\frac{1}{k}}
\left(
\sum_{\{i_1,i_2,\dots,i_k\}\in E(H)} u_{i_1}^k u_{i_2}^k \cdots u_{i_k}^k
\right)^{\frac{1}{k}}.
\]
Let $\mathbf{x}=(u_1^k,u_2^k,\dots,u_n^k)^\top $, then $\mathbf{x} \in \mathfrak{S}$ by $\sum\limits_{i=1}^{n} u_i^k = 1$.   Theorem \ref{th3.1} yields that
\[
\sum_{\{i_1,i_2,\dots,i_k\}\in E(H)} u_{i_1}^k u_{i_2}^k \cdots u_{i_k}^k
\le
\left(\frac{1}{\omega}\right)^k \binom{\omega}{k},
\]
then
\begin{align*}
\rho(\mathcal{A}) \le k m^{1-\frac{1}{k}}
\left(
\left(\frac{1}{\omega}\right)^k \binom{\omega}{k}
\right)^{\frac{1}{k}}.
\end{align*}
It follows from (\ref{eq4.2}) that
\begin{align*}
\left( \frac{1}{\omega} \right)^k \binom{\omega}{k}
\le \frac{1}{k!} \cdot \left( 1 - \frac{1}{\omega} \right)^{k-1}.
\end{align*}
Therefore,
\[
\omega \ge \frac{1}{1-\sqrt[k-1]{k!m^{1-k}\left(\frac{\rho(\mathcal{A})}{k}\right)^k}}.
\]
\end{proof}

To compare the two lower bounds given in Theorems~\ref{th4.1}
and~\ref{th4.3}, we consider the following two $3$-uniform
hypergraphs.

\begin{example}\label{ex4.4}
Let $H_1$ be the $3$-uniform hypergraph with
$
V(H_1)=\{1,2,3,4,5,6,7\}
$
and
$
E(H_1)=
\big\{
\{1,2,3\},\{1,2,4\},\{1,3,4\},\{2,3,4\},
\{5,6,7\}
\big\}.
$
It is easy to verify that
\[
\omega_2(H_1)=\omega(H_1)=4.
\]
Let  $
\rho(\mathcal{A}_{H_1})
$ denote the spectral radius of the adjacency tensor
$
\mathcal{A}_{H_1}
$ for $H_1$. By solving the corresponding eigenvalue equations, we obtain
$
\rho(\mathcal{A}_{H_1})=3,
$
and the  corresponding  nonnegative eigenvector
$
\mathbf{u}
=
\left(
4^{-\frac{1}{3}},4^{-\frac{1}{3}},4^{-\frac{1}{3}},4^{-\frac{1}{3}},0,0,0
\right)^\top,
$
$\|\mathbf{u}\|_3=1$.
Hence
\[
\|\mathbf{u}\|_1=4^{2/3}.
\]
By Theorem~\ref{th4.1}, we obtain
\[
\omega({H_1})
\ge
\frac{1}
{1-\sqrt{\dfrac{2\rho(\mathcal{A}_{H_1})}
{\|\mathbf{u}\|_1^3}}}
=
\frac{1}{1-\sqrt{3/8}}
\approx 2.5798.
\]
Since $\omega({H_1})$ is an integer, this implies
\[
\omega({H_1})\ge 3.
\]

On the other hand, since $m=5$, Theorem~\ref{th4.3} gives
\[
\omega({H_1})
\ge
\frac{1}
{1-\sqrt{
3!\cdot 5^{-2}
\left(\frac{\rho(\mathcal{A}_{H_1})}{3}\right)^3}}
=
\frac{1}{1-\sqrt{6/25}}
\approx 1.9604,
\]
which yields only
\[
\omega({H_1})\ge 2.
\]
Thus, for ${H_1}$, the bound in Theorem~\ref{th4.1} is sharper.

\vspace{3mm}
Next, let ${H_2}$ be the $3$-uniform hypergraph with
$
V({H_2})=\{1,2,3,4,5,6\}
$
and
$
E({H_2})=
\big\{
\{1,4,5\},\{1,4,6\},\{1,5,6\},\{4,5,6\},
\{1,2,5\},\{1,3,5\}
\big\}.
$
It is easy to verify that
\[
\omega_2({H_2})=\omega({H_2})=4.
\]
Let  $
\rho(\mathcal{A}_{H_2})
$ denote the spectral radius of adjacency tensor
$
\mathcal{A}_{H_2}
$ for $H_2$. By solving the corresponding eigenvalue equations, we obtain
$
\rho(\mathcal{A}_{H_2})\approx 3.5803
$,
and the corresponding nonnegative eigenvector
$
\mathbf{u}\approx
(0.6498,0.3434,0.3434,0.5699,0.6498,0.5699)^\top,
$
$\|\mathbf{u}\|_3=1$.
 Hence
\[
\|\mathbf{u}\|_1\approx3.1262.
\]

Theorem~\ref{th4.1} gives
$
\omega({H_2})\ge 1.9384,
$
and hence
\[
\omega({H_2})\ge 2.
\]
In contrast, since $m=6$, Theorem~\ref{th4.3} gives
$
\omega({H_2})\ge 2.1379,
$
and therefore
\[
\omega({H_2})\ge 3.
\]
Thus, for ${H_2}$, the bound in Theorem~\ref{th4.3} is sharper.

Consequently, the bound in Theorem~\ref{th4.1} is sharper for $H_1$,
whereas the bound in Theorem~\ref{th4.3} is sharper for $H_2$.
Thus, the two bounds are incomparable in general.
\end{example}

By Theorem~\ref{th4.1}  and the proof of the  Theorem~\ref{th4.3}, we directly get the following corollary.

\begin{cor}
Let $H$ be a $k$-uniform  hypergraph with $n$ vertices, $m$ edges and  clique number $\omega$. If $\omega_2(H)=\omega$,
then
\[
\rho(\mathcal A) \le \min\big\{k\|\mathbf{u}\|_1^k\left(\frac{1}{\omega}\right)^k \binom{\omega}{k},
~ k m^{1-\frac{1}{k}}
\left(
\left(\frac{1}{\omega}\right)^k \binom{\omega}{k}
\right)^{\frac{1}{k}}\big\},
\]
where $\mathbf{u}$ is a nonnegative eigenvector corresponding to $\rho(\mathcal{A})$ and $\|\mathbf{u}\|_k=1$.
\end{cor}

\subsection{Bounds based on the signless Laplacian and Laplacian tensors}

By Lemma \ref{perron-Frobenius}, the spectral radius
$\rho(\mathcal{Q})$ of the  signless Laplacian tensor $\mathcal{Q}$  is nonnegative and  there is  a nonnegative eigenvector corresponding to it. Let $\mathbf{v}=(v_1,v_2,\dots,v_n)^\top $ be a nonnegative eigenvector corresponding to $\rho(\mathcal{Q})$ and
  $\|\mathbf{v}\|_k= (\sum\limits_{i=1}^{n} v_i^k)^{\frac{1}{k}}=1 $.
Then we have the following  lower bound for the clique number of a
$k$-uniform hypergraph in terms of  $\rho(\mathcal{Q})$.
\begin{thm}\label{th4.5}
Let $H$ be a $k$-uniform hypergraph with clique number $\omega$ and maximum degree $\Delta$. If $\omega_2(H)=\omega$,
then
\begin{align*}
   \omega \ge \frac{1}{1 - \sqrt[k-1]{\dfrac{(k-1)!(\rho(\mathcal{Q})-\Delta)}{\|\mathbf{v}\|_1^k}}},
\end{align*}
where $\mathbf{v}$ is a nonnegative eigenvector corresponding to $\rho(\mathcal{Q})$ and $\|\mathbf{v}\|_k=1$.
\end{thm}
\begin{proof}
Let $\mathbf{y}=\frac{\textbf{v}}{\|\mathbf{v\|_1}}$, then  $\mathbf{y} \in \mathfrak{S}$. Thus, by Theorem \ref{th3.1},  we have
\begin{align*}
k\left(\frac{1}{\omega}\right)^k \binom{\omega}{k}&\ge k\sum_{\{i_1,\dots,i_k\}\in E(H)} y_{i_1}\cdots y_{i_k}
 \\
 &=  \sum_{i_1,\dots,i_k\in [n]} a_{i_1\dots i_k} y_{i_1}\cdots y_{i_k}       \\
 &= \sum_{i_1,\dots,i_k\in [n]} (a_{i_1\dots i_k}+d_{i_1\dots i_k})y_{i_1}\cdots y_{i_k}-\sum_{i_1,\dots,i_k\in [n]}d_{i_1\dots i_k}y_{i_1}\cdots y_{i_k}
 \\
&= \sum_{i_1,\dots,i_k\in [n]} (a_{i_1\dots i_k}+d_{i_1\dots i_k})\frac{v_{i_1}}{\|\mathbf{v}\|_1}\cdots
\frac{v_{i_k}}{\|\mathbf{v}\|_1}-\sum_{i_1,\dots,i_k\in [n]}d_{i_1 \dots i_k}\frac{v_{i_1}}{\|\mathbf{v}\|_1}\cdots
\frac{v_{i_k}}{\|\mathbf{v}\|_1} \\
&= \frac{1}{\|\mathbf{v}\|_1^k} \mathbf{v}^\top \mathcal{Q}\mathbf{v}^{k-1}-\frac{1}{\|\mathbf{v}\|_1^k} \sum_{i_1,\dots,i_k\in [n]}d_{i_1\dots i_k}{v_{i_1}}\cdots {v_{i_k}}\\
\end{align*}
\begin{align*}\label{eq4.3}
&= \frac{\rho(\mathcal{Q})}{\|\mathbf{v}\|_1^k}-\frac{1}{\|\mathbf{v}\|_1^k} (d_{1}{v_{1}^{k}}+\cdots+d_{n}{v_{n}^{k}})~~~~~~~~~~~~~~~~~~~~~~~~~~~\\
&\ge \frac{\rho(\mathcal{Q})}{\|\mathbf{v}\|_1^k}-\frac{\Delta}{\|\mathbf{v}\|_1^k} ({v_{1}^{k}}+\cdots+{v_{n}^{k}})\\
&= \frac{\rho(\mathcal{Q})-\Delta}{\|\mathbf{v}\|_1^k}. \tag{4.3}
\end{align*}
It follows from the inequality  (\ref{eq4.2}) that
$
\left( \frac{1}{ \omega} \right)^k \binom{ \omega}{k}
\le \frac{1}{k!} \left( 1 - \frac{1}{ \omega} \right)^{k-1}.
$
Then from inequality  (\ref{eq4.3}), we get
\begin{align*}
   \omega \ge \frac{1}{1 - \sqrt[k-1]{\dfrac{(k-1)!(\rho(\mathcal{Q})-\Delta)}{\|\mathbf{v}\|_1^k}}}.
\end{align*}

\end{proof}

By Lemmas \ref{QandL} and  \ref{Qi}, we get
$
\Delta \le\rho(\mathcal{L}) \le \rho(\mathcal{Q})
$.
So  Theorem \ref{th4.5} yields that
\begin{align*}
   \omega \ge \frac{1}{1 - \sqrt[k-1]{\dfrac{(k-1)!(\rho(\mathcal{Q})-\Delta)}{\|\mathbf{v}\|_1^k}}} \ge \frac{1}{1 - \sqrt[k-1]{\dfrac{(k-1)!(\rho(\mathcal{L})-\Delta)}{\|\mathbf{v}\|_1^k}}}.
\end{align*}
Then we have the following corollary.

\begin{cor}\label{th4.6}
Let $H$ be a $k$-uniform hypergraph with clique number $\omega$ and  maximum degree $\Delta$. If $\omega_2(H)=\omega$,
then
\begin{align*}
   \omega \ge \frac{1}{1 - \sqrt[k-1]{\dfrac{(k-1)!(\rho(\mathcal{L})-\Delta)}{\|\mathbf{v}\|_1^k}}},
\end{align*}
where $\mathbf{v}$ is a nonnegative eigenvector corresponding to $\rho(\mathcal{Q})$ and $\|\mathbf{v}\|_k=1$.
\end{cor}

Next, we give another lower bound for the clique number of a
$k$-uniform hypergraph in terms of the spectral radius of its
signless Laplacian tensor $\rho(\mathcal{Q})$.
\begin{thm}\label{th4.7}
Let $H$ be a $k$-uniform hypergraph with $n$ vertices, $m$ edges, maximum degree $\Delta$ and clique number $\omega$. If $\omega_2(H)=\omega$,
then
\[
\omega \ge \frac{1}{1-\sqrt[k-1]{k! m^{1-k}\left(\frac{\rho(\mathcal{Q})- \Delta}{k}\right)^k}}.
\]
\end{thm}
\begin{proof}
Let $\mathbf{v}=(v_1,v_2,\dots,v_n)^\top$ be a nonnegative eigenvector corresponding to $\rho(\mathcal{Q})$ and
  $ \sum\limits_{i=1}^{n} v_i^k=1 $. So
\begin{align*}
\rho(\mathcal{Q}) &= \mathbf{v}^\top (\mathcal{D}+\mathcal{A}) \mathbf{v}^{k-1} \\
&= \sum_{i=1}^{n} d_i v_i^k + k\sum_{\{i_1,i_2,\dots,i_k\}\in E(H)} v_{i_1}v_{i_2}\cdots v_{i_k} \\
&\le \Delta + k\sum_{\{i_1,i_2,\dots,i_k\}\in E(H)} v_{i_1}v_{i_2}\cdots v_{i_k}.
\end{align*}
Similar to the proof of Theorem \ref{th4.3}, we have
\begin{align*}\label{eq4.4}
\rho(\mathcal{Q}) \le \Delta + k m^{1-\frac{1}{k}}
\left(
\left(\frac{1}{ \omega}\right)^k \binom{ \omega}{k}
\right)^{\frac{1}{k}}.
\end{align*}
It follows from the inequality  (\ref{eq4.2}) that
$
\left( \frac{1}{ \omega} \right)^k \binom{ \omega}{k}
\le \frac{1}{k!} \left( 1 - \frac{1}{ \omega} \right)^{k-1}.
$
By Lemma \ref{QandL} and  \ref{Qi},  we have $\rho(\mathcal{Q}) - \Delta\geq0$.
Therefore,
\[
\omega \ge \frac{1}{1-\sqrt[k-1]{k!m^{1-k}\left(\frac{\rho(\mathcal{Q})- \Delta}{k}\right)^k}}.
\]
\end{proof}

By Lemmas \ref{QandL},  \ref{Qi} and Theorem \ref{th4.7},
  we have the following corollary.
\begin{cor}\label{cor4.8}
Let $H$ be a $k$-uniform hypergraph with $n$ vertices, $m$ edges and maximum degree $\Delta$, clique number $\omega$. If $\omega_2(H)=\omega$,
then
\[
\omega \ge \frac{1}{1-\sqrt[k-1]{k! m^{1-k}\left(\frac{\rho(\mathcal{L})- \Delta}{k}\right)^k}}.
\]
\end{cor}

\begin{rem}
When $k=2$, Corollary  \ref{cor4.8} coincides with  Theorem \ref{th1.3} for graphs. Therefore, Corollary  \ref{cor4.8} extends the result of graphs  to  $k$-uniform hypergraphs.
\end{rem}

By Theorem~\ref{th4.5}  and the proof of the  Theorem~\ref{th4.7}, we directly get the following corollary.
\begin{cor}
Let $H$ be a $k$-uniform hypergraph with $n$ vertices, $m$ edges, maximum degree $\Delta$ and clique number $\omega$. If $\omega_2(H)=\omega$,
then
\begin{align*}
\rho(\mathcal{L})&\le \rho(\mathcal Q) \le \min\big\{\Delta+k\|\mathbf{v}\|_1^k\left(\frac{1}{\omega}\right)^k \binom{\omega}{k},~~ \Delta + k m^{1-\frac{1}{k}}
\left(
\left(\frac{1}{ \omega}\right)^k \binom{ \omega}{k}
\right)^{\frac{1}{k}}\big\},
\end{align*}
where $\mathbf{v}$ is a nonnegative eigenvector corresponding to $\rho(\mathcal{Q})$ and $\|\mathbf{v}\|_k=1$.
\end{cor}

\section{Bounds under  an edge-number condition}

The proofs in Section~4 use the condition $\omega_2(H)=\omega$ only to ensure that
\[
L(H)=\left(\frac{1}{\omega}\right)^k\binom{\omega}{k}.
\]
Therefore, by Lemma~\ref{peng}, the same arguments can be applied
under the following edge-number condition, yielding analogous results.

\begin{cor}
Let $H$ be a $k$-uniform hypergraph with $n$ vertices, $m$ edges,
maximum degree $\Delta$, and clique number $\omega$. Suppose that
$
m\leq \binom{\omega}{k}+2(\omega+1-k).
$
Then
\begin{align*}
 \omega &\ge \frac{1}{1 - \sqrt[k-1]{\dfrac{(k-1)!\rho(\mathcal{A})}{\|\mathbf{u}\|_1^k}}},\\
 \omega &\ge \frac{1}{1-\sqrt[k-1]{k!m^{1-k}\left(\frac{\rho(\mathcal{A})}{k}\right)^k}},\\
 \omega &\ge \frac{1}{1 - \sqrt[k-1]{\dfrac{(k-1)!(\rho(\mathcal{Q})-\Delta)}{\|\mathbf{v}\|_1^k}}}\ge \frac{1}{1 - \sqrt[k-1]{\dfrac{(k-1)!(\rho(\mathcal{L})-\Delta)}{\|\mathbf{v}\|_1^k}}},\\
\omega &\ge \frac{1}{1-\sqrt[k-1]{k!\left(\frac{\rho(\mathcal{Q})- \Delta}{k}\right)^k m^{1-k}}} \ge \frac{1}{1-\sqrt[k-1]{k!\left(\frac{\rho(\mathcal{L})- \Delta}{k}\right)^k m^{1-k}}}.
\end{align*}
Moreover
\begin{align*}
 \rho(\mathcal A) &\le \min\big\{k\|\mathbf{u}\|_1^k\left(\frac{1}{\omega}\right)^k \binom{\omega}{k},
~ k m^{1-\frac{1}{k}}
\left(
\left(\frac{1}{\omega}\right)^k \binom{\omega}{k}
\right)^{\frac{1}{k}}\big\},\\
 \rho(\mathcal{L})&\le \rho(\mathcal Q) \le \min\big\{\Delta+k\|\mathbf{v}\|_1^k\left(\frac{1}{\omega}\right)^k \binom{\omega}{k},~~ \Delta + k m^{1-\frac{1}{k}}
\left(
\left(\frac{1}{ \omega}\right)^k \binom{ \omega}{k}
\right)^{\frac{1}{k}}\big\}.
\end{align*}
Here, $\mathbf{u}$ and $\mathbf{v}$ are nonnegative eigenvectors corresponding to $\rho(\mathcal{A})$ and $\rho(\mathcal{Q})$, respectively, with $\|\mathbf{u}\|_k=\|\mathbf{v}\|_k=1$.
\end{cor}

\begin{rem}
The results in Section~4 depend essentially on the identity
\[
L(H)=\left(\frac{1}{\omega}\right)^k\binom{\omega}{k}.
\]
Hence, whenever this identity holds under other conditions on the hypergraph, the same conclusions on clique number and spectral bounds can be obtained.
\end{rem}

\end{spacing}

\vspace{3mm}
\noindent
\textbf{
Acknowledgements.}
\vspace{3mm}

This work is supported by
Science Research Project of Hebei Education Department (No. QN2025023),
Young Innovative Talents Project of General Colleges and
Universities in Guangdong Province (No. 2022KQNCX018).


\vspace{3mm}
\noindent
\bibliographystyle{elsarticle-num}

\end{CJK*}
\end{document}